\documentclass[12pt]{amsart}
\usepackage{amsmath,amssymb,amsthm,amsfonts}

\usepackage{graphicx}
\usepackage{amssymb}
\usepackage{amsmath}
\usepackage{color}
\usepackage{times}

\usepackage[unicode,bookmarks,colorlinks]{hyperref}
\hypersetup{
    linkcolor=brickred,
}

\definecolor{mahogany}{cmyk}{0, 0.77, 0.87, 0}
\definecolor{salmon}{cmyk}{0, 0.53, 0.38, 0}
\definecolor{melon}{cmyk}{0, 0.46, 0.50, 0}
\definecolor{yellowgreen}{cmyk}{0.44, 0, 0.74, 0}
\definecolor{brickred}{cmyk}{0, 0.89, 0.94, 0.28}
\definecolor{OliveGreen}{cmyk}{0.64, 0, 0.95, 0.40}
\definecolor{RawSienna}{cmyk}{0, 0.72, 1.0, 0.45}
\definecolor{ZurichRed}{rgb}{1, 0, 0} 

\usepackage{amsmath,amstext,amssymb,amsopn,amsthm}
\usepackage{amsmath,amssymb,amsthm}
\usepackage[mathscr]{eucal}

\begin{document}

\newtheorem{lemma}[thm]{Lemma}
\newtheorem{remark}{Remark}
\newtheorem{proposition}{Proposition}
\newtheorem{theorem}{Theorem}[section]
\newtheorem{deff}[thm]{Definition}
\newtheorem{case}[thm]{Case}
\newtheorem{prop}[thm]{Proposition}
\numberwithin{equation}{section}

\numberwithin{remark}{section}
\numberwithin{proposition}{section}

\newcommand{\gap}{\lambda_{2,D}^V-\lambda_{1,D}^V}
\newcommand{\gapR}{\lambda_{2,R}-\lambda_{1,R}}
\newcommand{\bD}{\mathrm{I\! D\!}}
\newcommand{\calD}{\mathcal{D}}
\newcommand{\calA}{\mathcal{A}}

\newcommand{\conjugate}[1]{\overline{#1}}
\newcommand{\abs}[1]{\left| #1 \right|}
\newcommand{\cl}[1]{\overline{#1}}
\newcommand{\expr}[1]{\left( #1 \right)}
\newcommand{\set}[1]{\left\{ #1 \right\}}

\newcommand{\calC}{\mathcal{C}}
\newcommand{\calE}{\mathcal{E}}
\newcommand{\calF}{\mathcal{F}}
\newcommand{\Rd}{\mathbb{R}^d}
\newcommand{\BR}{\mathcal{B}(\Rd)}
\newcommand{\R}{\mathbb{R}}
\newcommand{\bH}{\mathbb{H}}

\newcommand{\al}{\alpha}
\newcommand{\RR}[1]{\mathbb{#1}}
\newcommand{\bR}{\mathrm{I\! R\!}}
\newcommand{\ga}{\gamma}
\newcommand{\om}{\omega}
\newcommand{\A}{\mathbb{A}}

\newcommand{\bb}[1]{\mathbb{#1}}
\newcommand{\bI}{\bb{I}}
\newcommand{\bN}{\bb{N}}

\newcommand{\uS}{\mathbb{S}}
\newcommand{\M}{{\mathcal{M}}}
\newcommand{\W}{{\mathcal{W}}}

\newcommand{\m}{{\mathcal{m}}}

\newcommand {\mac}[1] { \mathbb{#1} }

\newcommand{\bC}{\Bbb C}
\date{December 21, 2010}

\title[Martingales]{Donald Burkholder's work in martingales and analysis}\thanks{To appear in ``Collected works of D.L. Burkholder"}
\author{Rodrigo Ba\~nuelos}
\thanks{R. Ba\~nuelos is supported in part  by NSF Grant
\# 1005844-DMS}
\address{Department of Mathematics, Purdue University, West Lafayette, IN 47907}
\email{banuelos@math.purdue.edu}
\author{Burgess Davis}
\address{Department of Mathematics and Department of Statistics,  Purdue University, West Lafayette, IN 47907}
\email{bdavis@stat.purdue.edu}
\maketitle
\tableofcontents
\section{Introduction}
The two mathematicians who have most advanced martingale theory in the last seventy years are Joseph Doob and Donald Burkholder. Martingales as a remarkably flexible tool are used throughout probability and its applications to other areas of mathematics. They are central to modern stochastic analysis. And martingales, which can be defined in terms of fair games, lie at the core of mathematical finance. Burkholder's research has profoundly advanced not only martingale theory but also, via martingale connections, harmonic and functional analysis.

The work of Burkholder and Gundy on martingales in the late sixties and early seventies, which followed Burkholder's seminal 1966 paper { \it  Martingale Transforms} \cite{Bur11}, led to applications in analysis which revolutionized parts of this subject. Burkholder's outstanding work in the geometry of Banach spaces, described by Gilles Pisier in this volume, arose from his extension of martingale inequalities to settings beyond Hilbert spaces where the square function approach used in \cite{Bur11} fails. His work in the eighties and nineties on martingale inequalities with emphasis on identifying best constants has become of great importance recently in the investigations of two well known open problems. One of these concerns optimal $L^p$ bounds for a singular integral operator (the two dimensional Hilbert transform) and their ramifications in quasiconformal mappings. The other relates to a longstanding conjecture in the calculus of variations dealing with rank-one convex and quasiconvex functions. These conjectures, which have received much attention in recent years largely  due to the beautiful and original techniques developed by Burkholder in his work on sharp martingale inequalities, come from fields which on the surface are far removed from martingales. 

We will describe in some detail a remarkable technique discovered by Burkholder and Gundy, which shows how certain integral inequalities between two nonnegative functions on a measure space follow from inequalities involving only parts of their distribution.  This seemingly simple but incredibly elegant technique, often, and here, referred to  as { ``the good--$\lambda$ method,"} revolutionized the way probabilists and analysts think of norm comparison problems.  It is now widely used in areas of mathematics which involve integrals and operators.

It is interesting to note that since 1973, Burkholder has written only two papers with a co-author and that
he has written more than one paper only with Richard Gundy.
The papers \cite{Bur14}  of Burkholder and Gundy and \cite{Bur16} of Burkholder, Gundy, and Silverstein are exceptionally important.  The results of \cite{Bur14} include the good--$\lambda$ inequalities and fundamental integral inequalities comparing the maximal function and the square function, or quadratic variation, of martingales having controlled jumps or continuous paths.
A very large share of the extensive applications of these kinds of martingale inequalities, both in probability and other areas of mathematics, involve continuous path martingales.
The paper \cite{Bur16} strikingly improved and completed work of Hardy and Littlewood on the characterization of the Hardy $H^p$ spaces via the integrability of certain maximal functions. While probabilistic techniques had already gained the respect of many analysts studying harmonic functions and potential theory, due in part to earlier work of Doob, Kakutani, Wiener and others, this landmark paper had a profound influence in harmonic analysis and propelled many analysts to learn probability.

      The next section begins with a brief introduction to the good--$\lambda$ method, in the context of its original application to martingales. We then trace the rest of the development of the theory of martingale square functions and transforms in the late sixties and early seventies, pioneered by Burkholder. We follow this with a discussion of \cite{Bur16} and the subsequent  study of $H^p$ theory by a number of researchers, and much more on the surprisingly rich good-$\lambda$ inequalities. In the final sections we discuss Burkholder's later work on sharp martingale inequalities and some of the remarkable spread of his ideas over other areas of mathematics.

\section{Martingale Inequalities}

Brownian motion stopped at a stopping time is a continuous martingale, and the continuous martingale inequalities of \cite{Bur14} follow from their validity just for stopped Brownian motions.
We will elaborate on this later.
We will  use  $B=\{B_t; \ t\geq 0\}$ to denote standard Brownian motion.  This means that $B$ is a stochastic process with continuous paths, that its increment $B_t-B_s$ over the interval $[s,t]$ has a normal distribution with mean 0 and variance $t-s$, that its increments over each of a collection of disjoint intervals are independent, and that $B_0=0$.
We recall that if the random variable $\tau$ is a stopping time for $B$ then $\tau\geq 0$ and if $P(\tau>s)>0$ and $t>s$,
the conditional distribution, given $\tau>s$, of $B_t-B_s$ is normal with mean zero and variance $t-s$.
The maximal function of $B$ up to the stopping time $\tau$ will be denoted by $B_{\tau}^* =\sup\{|B_s|\colon 0\leq s<\tau\}$.  The following theorem is from \cite{Bur14}.

\begin{theorem}\label{thm1}
\ Let $\Phi$ be a continuous nondecreasing function on $[0,\infty)$ satisfying  $\Phi(0)=0$ and $\Phi(2\lambda)\leq K\Phi(\lambda)$,  $\lambda\geq 0$, for some constant $K$.
Then there are positive constants $c$ and $C$, which depend only on $K$, such that for any stopping time $\tau$ for B,
\begin{equation}\label{norm1}
cE\Phi(\sqrt{\tau})\leq E\Phi(B^*_{\tau}) \leq CE\Phi(\sqrt{\tau}).
\end{equation}
\end{theorem}

\begin{remark} Two important examples of functions $\Phi$ satisfying this ``moderate" growth property are $\Phi(x)=x^p$, $0<p<\infty$, and $\Phi(x)=x+x\ln^{+}(x)$.
\end{remark}

To illustrate the good-$\lambda$ method used by Burkholder and Gundy in \cite{Bur14}, we give a direct proof of the left hand side of (\ref{norm1}) in the case $\Phi(x)=x$ which gives $c=\frac{1}{1200}$.  This proof, which requires virtually no specialized knowledge, is a slight alteration of the proof in \cite{Bur14}, as it uses summation rather than integration. Later,  in Theorems \ref{thm2}and \ref{thm3}, we present a general form of the good--$\lambda$ method, together with inequalities for stopped Brownian motion, which imply Theorem \ref{thm1}.

Denote the integers by ${\bf Z}$.
Let $a_k\geq 0$, $k\in \mathbf Z,$ satisfy $\lim_{k\to -\infty}a_{k}=0$ and $a_{k+1}\leq 2a{_k}$.
For  $0<r<1$, let $J(r)=\{ k\colon a_{k+1}>ra_k\}$.
If $k$ is in $J(r)$, but none of $k+i$, for $1\leq i\leq m$, are in $J(r)$, then
$$\sum_{i=1}^ma_{k+i} \leq a_{k+1}(1+r+r^2+\ldots r^{m-1}) <\frac{2a_k}{1-r},$$ which implies
\begin{equation}\label{tag2} 
\sum_{k\in J(r)} a_k\geq\frac{1-r}{3-r}\sum_{k\in \mathbf Z}a_k.
\end{equation}
The $k$ in $J(r)$ are the ``good'' $k$.
Now nonnegative random variables $X$ satisfy
\begin{equation}\label{tag3}
EX\leq\sum_{k\in \mathbf Z} 2^k P(X\geq 2^k)\leq 2EX.
\end{equation}
If $N$ is a standard normal random variable then, using tables or that the density of $N$ is bounded by $\frac{1}{\sqrt{2\pi}}$, we get $P\left(|N|<\frac{1}{10}\right) <\frac{1}{12}$, so for an
event $A$,
\begin{equation}\label{tag4}
P\left(|N|\geq\frac{1}{10}, A\right)\geq\frac{1}{6}, \text{ if }  P(A)\geq\frac{1}{4}.
\end{equation}
Let $A_k=\{\sqrt{\tau}\geq 2^k\}, k\in \mathbf Z$, and let $J=\{ k\colon P(A_{k+1})\geq\frac{P(A_k)}{4} \}$.
The left hand side of (\ref{tag3}), and (\ref{tag2}) with $r=1/2$ and $a_k=2^k P(A_k)$, give
\begin{equation}\label{tag5}
E\sqrt{\tau}\leq 5\Sigma_{k\in J} 2^k P(A_k).
\end{equation}
Since $2B^*_t\geq |B_a-B_b|$, if $0\leq a\leq b\leq t$,  $2B^*_{\tau} \geq | B_{2^{2(k+1)}}-B_{2^{2k}}|$ on $A_{k+1}$.
With (\ref{tag4}) this gives
$$
P(20B_{\tau}^*\geq2^k)\geq P(2B_{\tau}^* \geq \frac{1}{10} 2^k\sqrt{3}, A_{k+1})\geq \frac{1}{6} P(A_k),\  k\in J,
$$
which with the right side of (\ref{tag3}) and (\ref{tag5}) yields
$$
2E20B_{\tau}^*\geq \sum_{k\in \mathbf Z} P(20B^*_{\tau}\geq 2^k) 2^k
\geq\frac{1}{6}\sum_{k\in J} P(A_k) 2^k\geq \frac{1}{30} E\sqrt{\tau}.
$$

As noted in \cite{Bur20}, Skorohod and others had before \cite{Bur14} proved the inequalities (\ref{norm1}) for the case $\Phi(x)=x^p$, $p\geq 2$, and P. W.~Millar \cite{Mil},  using results of \cite{Bur11}, extended these to all $p>1$.
Also A. A.~Novikov \cite{Nov}, working independently of \cite{Bur14}, used stochastic calculus to study questions raised by Millar's paper and proved some interesting results related to those of \cite{Bur14}.

The growth condition on $\Phi$ involving $K$ of Theorem \ref{thm1} is necessary for the truth of any of the inequalities in (\ref{norm1}), in the sense that if $\Phi$ is a continuous nondecreasing function which does not satisfy this condition for any $K$ there are stopping times $\tau$ for $B$ such that (either) one of $E\Phi(B_\tau^*)$, $E\Phi(\sqrt{\tau})$ is finite and the other is infinite.

Next we turn to discrete time martingales.
After a very brief history of martingales before Doob we provide an overview of the work of the late sixties and seventies involving the martingale square function.
More general results, with proofs and extensive references, may be found in Burkholder's Wald Memorial Lecture paper \cite{Bur20}.
We have tried to be true to the spirit if not the letter of the papers we describe.

Paul L\'evy defined martingales without the name, which was given by Doob.
Before martingales were formally defined, several probabilists other than L\'evy, and several analysts, worked on objects that were martingales.
For example, R.E.A.C. Paley \cite{Pal} proved an inequality for the Haar system which is a special case of the results of Burkholder in his 1966 paper.  (See \cite{Bur40} for a sharp version of the Paley result.)
Although the definition of martingales was made by a probabilist, there is no reason it couldn't have come from an analyst instead.
Sequences of piecewise constant functions on the Lebesgue unit interval which are martingales seem now a natural generalization of Haar series, and are in a distributional sense all of the discrete (as described in the next paragraph) martingales.
Of course, there's nothing like hindsight to clarify thinking.
In another direction Courant, Fredricks, and Lewy in 1928 \cite{CouFreLew} used ideas related to martingale ideas, although without randomness, to study harmonic functions, in the paper  which introduced the finite element method for numerical approximation of solutions of partial differential equations.

We begin with a description of martingales when time is discrete and the random variables which compose them are discrete, that is, have a discrete distribution.
A sequence of discrete random variables $\{D_i, i\geq 0\}$, is a martingale difference sequence if each $D_i$  has  finite expectation and if for $n>0$, 
\begin{equation}\label{desmart}
E(D_n|D_i=a_i, 0\leq i<n)=0,  \text{ if } P(D_i=a_i,0\leq i<n)>0.
\end{equation}
We may think of a gambler as having initial stake $D_0$ and playing a sequence of games with the amount won or lost upon playing the $i$th game being $D_i$.
The game the gambler plays at time $n$ may differ depending on her initial stake and her history of wins and losses in the first $(n-1)$ games, but the game she will play at time $n$ always has expectation zero.
Of course the sequence of partial sums $\sum_{i=0}^n D_i,\ n\geq 0$, is the martingale corresponding to the difference sequence.
Except for some technicalities, the study of these fully discrete martingales is invested with all the difficulties connected with the more general discrete time martingales described below.

We start with a probability space $(\Omega, \mathcal{A}, P)$, and a sequence $\mathcal{A}_0, \mathcal{A}_1, \dots$ of $\sigma$-fields contained in $\calA$ such that $\calA_n\subset \calA_{n+1},\ n\geq 0$.
A sequence of random variables $f=\{f_n, n\geq 0\}$, on $\Omega$ is a martingale with respect to these $\sigma$-fields if each $f_i$ is $\calA_i$ measurable and integrable and if $E(d_i|\calA_{i-1})=0,\ i>0$, where $d_i=f_i-f_{i-1}$, $d_0=f_0$.

The maximal function $f^*$ and the square function $S(f)$ of a martingale $f$ with difference sequence $d$ are defined by
\begin{equation*}
f^*=\sup_n |f_n|, \hskip.3cm \text{and} \hskip.3cm  S(f)=\left(\sum_{i=0}^\infty d_i^2\right)^{1/2}.
\end{equation*}
The conditional version of Jensen's inequality implies that  $E|f_n|^p$ is nondecreasing for $1\leq p<\infty$, so that
$\lim_{n\rightarrow\infty}(E|f_n|^p)^{1/p}$ makes sense.
This limit is called the $L^p$ norm of $ f$ and denoted $\|f\|_p$ and if $\|f\|_p<\infty$, $f$ is said to be $L^p$ bounded. 
The celebrated maximal inequalities of Doob assert that  for any martingale $ f$,
\begin{equation}\label{doobweak}
P\{f^*>\lambda\}\leq \frac{1}{\lambda} \|f\|_1, \,\,\, \lambda>0, 
\end{equation}
and 
\begin{equation}\label{doobstrong}
\|f^*\|_p \leq \frac{p}{p-1} \|f\|_p, \,\,\, 1<p<\infty.
\end{equation}

The martingale differences $\{d_n, n\geq 0\}$ of any $L^2$ bounded martingale ${ f}$ are orthogonal and it follows trivially that  $\|S(f)\|_2=\|f\|_2$,  which implies $S(f)<\infty$ a.s. (almost surely).
D. Austin in his 1966 paper \cite{Aus} strengthened this by showing that $S(f)<\infty$ a.s.~if $f$ is an $L^1$ bounded martingale.
Burkholder in \cite{Bur11} in turn strengthened Austin's result by proving that the operator $f\to S(f)$ is weak-type $(1,1)$.  That is, he proved that there is a universal constant $C$  such that if $ f$ is an $L^1$ bounded martingale then
\begin{equation}\label{weakS}
P\{S(f)>\lambda\}\leq \frac{C}{\lambda}\|f\|_1, \,\,\, \lambda>0.
\end{equation}

That (\ref{weakS}) might be true was informed by both Austin's result and earlier work of Burkholder himself, especially his paper \cite{Bur10}, {\it  Maximal inequalities as necessary conditions for almost everywhere convergence}.  We sketch just the start of this argument.
If (\ref{weakS}) does not hold for any $C$ then neither does it hold for any $C$ for all martingales having a finite index (time) set and with initial value (i.e.~$d_0$) equal to 0.
So there is a sequence ${g_1, g_2, \ldots}$ of such martingales satisfying
$$P\{S(g_i)\geq y_n\}\geq\frac{C_n}{y_n} \|g_i\|_1,$$ for some positive numbers $y_n$ and $C_n$ such that $C_n\rightarrow\infty$, as $n\rightarrow\infty$.
However, from these martingales a martingale which has finite $L^1$ norm and almost surely infinite square function can be constructed, by putting independent copies of the martingales $k_{n_i} g_{n_i}$ sequentially, where $n_i, i\geq 0$, is a sequence of positive integers (with some integers repeated) and the $k_i$ are constants. This would contradict Austin's theorem.
Burkholder in \cite{Bur20} gives an elementary  proof of (\ref{weakS}) with $C=3$.  Later, in \cite{Bur31}, he proves the inequality with
$C=2$.
In  \cite{Cox} Cox proves the inequality with $C=e^{1/2}$ and shows that this is best possible.

Burkholder also showed in  \cite{Bur11} that there are constants $c_p$ and $C_p$ such that for any martingale
\begin{equation}\label{strongS}
c_p \|f\|_p\leq \|S(f)\|_p \leq C_p\|f\|_p, \,\,\, 1<p<\infty.
\end{equation}
This starts to suggest  the interchangeability of $f^*$ and $S(f)$ that appears in so many of the results of \cite{Bur14}.
An example of such interchangeability, which follows quickly from results of \cite{Bur11}, comes from combining (\ref{doobstrong}) and (\ref{strongS}) to get the existence of positive constants $k_p$ and $K_p$ such that for any martingale, 
\begin{equation}\label{strongmaxS}
k_p \|f^*\|_p\leq \|S(f)\|_p\leq K_p\|f^*\|_p , \,\,\,  1<p<\infty.
\end{equation}

If $f$ and $g$ are martingales with difference sequences $d$ and $e$, respectively, satisfying $|e_i|\leq |d_i|$ for all $i$, we say that  $g$ is differentially subordinate to $ f$.
In this case $S(g)\leq S(f)$ and so (\ref{strongS}) gives that if $f$ is in $L^p$ for some $p>1$, then so is $g$.   Differential subordination sometimes occurs when $g$ and $f$ are both constructed from the same process or when one is constructed from the other, for example when $g$ is a {martingale transform} of $f$.
Martingale transforms are  discrete versions of It\^o integrals.
If $f$ is a martingale with respect to the sequence of $\sigma$--fields $\{\calA_i, i\geq 0\}$,
with difference sequence $d$, and if ${v}=\{v_i, i\geq 0\}$ is a sequence of random variables with the property that  $v_i$ is $\calA_{i-1}$  measurable, $i>0$ and $v_0$ is a constant (such a sequence of random variables is said to be predictable relative to the family of $\sigma$-fields $\{\calA_i\}$), then the process with difference sequence $v_{i}d_{i}$ is called a {martingale transform} of f and denoted $v \ast f$.
(In the completely discrete case of the gambler described above, if we take the $\sigma$--fields to be those generated by the random variables, then $v_n$ will be constant on each of the events  $\{D_i=a_i,\ 0\leq i<n\}$ conditioned on in (\ref{desmart}).
So the transformed gambler still always plays a fair game, derived from the original game by changing the stakes and/or the gambler changing places with her opponent.)
This  transformation of the martingale difference sequence  of $f$ may not yield a martingale because $v_{i}d_{i}$ may not be integrable. However,  under the assumption that  $|v_i|\leq1$  for all $i$, this problem does not arise, and the difference sequence  $\{v_id_i, i\geq 0\}$ generates a new martingale.
Furthermore $S(v\ast f)\leq S(f)$, so  (\ref{strongS}) implies that the operation $f\to v\ast f$ is bounded in $L^p$ for $1<p<\infty$.
The following theorem is from \cite{Bur11}.

\begin{theorem}  Let $f$ be a martingale on the sequence of $\sigma$-fields $\{\calA_i, i\geq 0\}$ with difference sequence $d$.  Let $\{v_i, i\geq 0\}$ be a predictable sequence with $|v_i|\leq 1$ a.s. for all $i$.  Then there are constants $C_1$ and $C_p$ such that 
\begin{equation}
P\{{(v\ast f)}^*>\lambda\}\leq \frac{C_1}{\lambda}\|f\|_1,\,\,\, \lambda >0
\end{equation}
and
\begin{equation}
\|v\ast f\|_p\leq C_p\|f\|_p, \,\,\, 1<p<\infty.
\end{equation}
\end{theorem}

A process $Y=\{Y_t,\ t\geq 0\}$ which has continuous paths is a martingale if for every positive integer $n$, the discrete time process $Y_{ n}=\{ Y_{\frac{k}{n}}, k\geq 0\}$, is a martingale.
The maximal function of $ Y$ equals $\sup_t|Y_t|$, while the analog of the square function, called the quadratic variation of $Y$ and denoted  $\langle Y\rangle$, can be defined as either the limit of the square functions of the $Y_{ n}$ for $n$ increasing sufficiently fast, or as the quantity that needs to be subtracted from the submartingale $Y^2$ to make it a martingale.   $\langle Y\rangle$ is also the stopping time that results when
$ Y$ is time changed to a stopped Brownian motion.
See \cite{RevYor}. 
The standard It\^o  integral with respect to Brownian
motion is a continuous martingale.  It is  interesting that all continuous martingales on the $\sigma$--fields generated by a
Brownian motion can be represented as It\^o integrals.

As mentioned earlier, it was proved in \cite{Bur14}  that if $\Phi$ is as in Theorem \ref{thm1} then for any continuous martingale $ Y$ which starts at $0$,
\begin{equation}\label{norm2}
c E\Phi\left(\langle Y\rangle^{1/2}\right)\leq E\Phi\left(Y^*\right)\leq
CE\Phi\left(\langle Y\rangle^{1/2}\right)
\end{equation}
where the constants $c$ and $C$ are the same as in (\ref{norm1}).
A number of integral inequalities for discrete martingales are proved in \cite{Bur14} using the good--$\lambda$ method and other techniques.
We confine ourselves here to those relating integrals of the maximal and the square function, but there are numerous other operators considered in \cite{Bur14}.
The next two theorems from \cite{Bur14} give versions of (\ref{norm2}) for large classes of martingales for which the jump size, that is the distributions of the $d_i$, is controlled.
The martingales to which the next theorem applies include stopped random walks where the jumps are iid with mean zero and finite variance.
In the following two theorems $f$ is a martingale with difference sequence $d$ with respect to the sigma fields $\calA_i$ as above.

\begin{theorem}
 Suppose
$E(d_i^2|\calA_{i-1})\leq M[E(|d_i | |\calA_{i-1})]^2$, $i>0$, for a positive constant $M$.
Then if $\Phi$ and $K$ are as Theorem 2.1 there exist constants $c$ and $C$ which depend only on $K$ and $M$ such that
\begin{equation}\label{norm3}
cE\Phi(f^*)\leq E\Phi((S(f))\leq CE\Phi(f^*).
\end{equation}
\end{theorem}

Another result from \cite{Bur14} is the following.
\begin{theorem} Suppose there is a number $M$ such that $|d_i|\leq M$ , $i\geq 0$, and that
$\Phi$ and $K$ are as in Theorem \ref{thm1}.
Then there is a constant $c$ depending only on $K$ such that
\begin{equation*}
E\Phi(S(f))  \leq  cE\Phi(f^*) +c\Phi(M)
\end{equation*}
and
\begin{equation*}
E\Phi(f^*)  \leq  cE\Phi(S(f)+c\Phi(M).
\end{equation*}
\end{theorem}
\begin{remark}
In particular, these inequalities imply that for martingales with uniformly bounded difference sequences if one of $E\Phi(f^*)$ or $E\Phi(S(f))$ is finite then the other is, which is often enough in applications.
\end{remark}

Old examples of Marcenkeiwicz and Zygmund, noted in \cite{Bur14}, show that for every $p\in(0,1)$, there are martingales for which  $\|f^*\|_p$ is finite but $\|S(f)\|_p$ is infinite.  Similar examples are given showing that for every $p\in(0,1)$ there are martingales for which $\|S(f)\|_p$ is finite but $\|f^*\|_p$ is infinite.
These examples show that neither side of (\ref{strongmaxS}) extends to $0<p<1$.

In \cite{Dav1}, (\ref{strongmaxS}) was extended to $p=1$ using a decomposition which trims the big jumps from a martingale, leaving a martingale with controllable jumps, which was handled with a method from \cite{Bur14}. See Garsia \cite{Garsia1} for a very different proof. 
Later, Burkholder, Davis and Gundy \cite{Bur18} used this decomposition and techniques related to the good--$\lambda$ method to prove that for every function  $\Phi$ as in Theorem \ref{thm1} which is also convex, there are constants $c$ and $C$, which depend only on $K$, such that for every martingale,
\begin{equation}
cE\Phi(f^*)\leq E\Phi(S(f))\leq CE\Phi(f^*).
\end{equation}
See Garsia \cite{Garsia2} and \cite{Garsia3} for a different approach to some of the inequalities of \cite{Bur18} and for related martingale inequalities.

\section{Martingale Inequalities and Hardy spaces}

We now turn to applications of (\ref{norm2}) and the good--$\lambda$ method to analysis, beginning by quoting the first seven lines of the Burkholder, Gundy and Silverstein paper \cite{Bur16} {\it  ``A maximal function characterization of the class $H^p.$''}

``Hardy and Littlewood have shown [\cite{HL}; also page 278 of vol. 1 of Zygmund's book \cite{Zyg}]  that if $F(z)$ is analytic in the unit disc $D=\{z\in \bC: |z|<1\}$, and if $\Omega_{\sigma}(\theta)$ is the Stoltz domain given by the interior of the smallest convex set containing the disc $\{z\in \bC: |z|<\sigma\}$ and the point $e^{i\theta}$, then
\begin{equation}\label{BGS1}
\int_0^{2\pi} \sup_{z\in\Omega_{\sigma}(\theta)} |ReF(z)|^p\,d\theta\leq
C_{\sigma, p} \sup_{0<r<1} \int_0^{2\pi}|F(re^{i\theta})|^pd\theta
\end{equation}
for all $p>0$, $0<\sigma <1$.
In this paper we prove the converse inequality which, together with the above theorem of Hardy and Littlewood, gives a maximal function characterization of the Hardy class $H^p$.''

The discussion of \cite{Bur16} we now provide is essentially a partial summary of Burkholder's survey paper \cite{Bur34}, included in this volume. 

We let $G(z) =u(z)+iv(z)$, where $G$ is continuous in $|z|\leq 1$ and analytic in $|z|<1$. 
We will refer to $u$ and $v$ as ``conjugate" harmonic functions.  The boundary function of $u$ is denoted by $U$, where
$U(e^{i\theta}) =u(e^{i\theta})$, with a similar notation for $V$.  We assume that $u(0)=v(0)=0$. Under this assumption $u$ and $v$ determine each other, but while this is easy to show, it is difficult to say much about many aspects of the connection between $u$ and $v$.  It is, however, easy to show that the (squares of the) $L^2$ norms of $U$ and $V$  are the same, that is   $\int_0^{2\pi}|U(e^{i\theta})|^2\,d\theta=\int_0^{2\pi}|V(e^{i\theta})|^2\,d\theta$. While equality need not hold if the exponent 2 is replaced by any other positive number $p$, M. Riesz showed that for $p>1$ the $L^p$ norm of $V$ cannot be more than a constant $C_p$ times the $L^p$ norm of $U$.  For $p=1$ even this does not hold, but Kolmogorov's weak type inequality provides
 some control of $V$ by the $L^1$ norm of $U$.
 
 Before we state the Burkholder-Gundy-Silverstein result from which the converse  of (\ref{BGS1}) is derived, we recall that for any function $f$ in the unit disc, the function 
\begin{equation*}\label{nontandisc}
N_{\sigma}f(\theta)=\sup\{|f(z)|:\,z\in \Omega_{\sigma}(\theta)\}
\end{equation*}
is called the nontangential maximal function of $f$.

The following theorem, proved in  \cite{Bur16}, immediately implies the appropriate, that is, with F(0)=0, 
converse of (\ref{BGS1}) asserted in the statement quoted above from \cite{Bur16}.

\begin{theorem}\label{BGS2}  Let $u$ and $v$ be conjugate harmonic functions with $u(0)=v(0)=0$. Let $\Phi$ and $K$ be as in Theorem \ref{thm1}.  
There  are positive constants  $c_{K, \sigma}$ and $C_{K, \sigma}$  depending only $K$  and $\sigma$ such that

\begin{equation}\label{BGS3}
c_{K, \sigma}\int_0^{2\pi} \Phi\left(N_{\sigma}u(\theta)\right) d\theta
\leq \int_0^{2\pi} \Phi\left(N_{\sigma}v(\theta)\right) d\theta
 \leq
C_{K, \sigma}\int_0^{2\pi}  \Phi\left(N_{\sigma}u(\theta)\right)d\theta.
\end{equation}
\end{theorem}

As pointed out by Burkholder in \cite{Bur34},  these inequalities in the cases $\Phi(x)=x^p, p>0,$ may be viewed as an extension of the Riesz inequalities to $0<p<\infty$. 

 The proof of this theorem uses (\ref{norm2}) applied to martingales obtained by composing harmonic functions with Brownian motion.  More precisely, we recall that  if $B$ is Brownian motion in $D$ and $\tau$ denotes the first time $B$ hits the unit circle, the process $u(B_{t\wedge \tau})$ is a martingale whenever $u$ is harmonic in the open unit disc and continuous in its closure.  The maximal function of this martingale, denoted by $u^*$, is called the { Brownian maximal function} of $u$.  If we continue with the normalization $u(0)=0$, then by the It\^o formula, the quadratic variation of this martingale is 
 \begin{equation}\label{brownsquare}
 \langle u\rangle_{t}=\int_{0}^{t\wedge \tau} |\nabla(B_s)|^2 \,ds. 
 \end{equation}
 This quantity is called the { Brownian square function} of $u$.  A key observation is that whenever  $u$ and $v$ are conjugate harmonic functions as above, then by the Cauchy--Riemann equations $|\nabla u|=|\nabla v|$ and hence $u$ and $v$ have the same Brownian square functions.  Thus the following theorem is an immediate consequence of (\ref{norm2}).

 \begin{theorem}\label{BGS4} Let $u$ and $v$ be conjugate harmonic functions. 
 Let $\Phi$ be as in the statement of Theorem \ref{thm1}.  There are constants $c$ and $C$ depending only on $K$ such that 
 \begin{equation}\label{BGS5}
 cE\Phi(u^*)\leq E\Phi(v^*)\leq CE\Phi(u^*).
 \end{equation}
 \end{theorem}
 
This is a {Brownian} version of (\ref{BGS3}).  
This theorem together with the following very surprising and beautiful result from \cite{Bur16} implies Theorem \ref{BGS2}. 

\begin{theorem}\label{BGS6} Let $u$ be a harmonic function in the unit disc. There are constants $c_{\sigma}$ and  $C_{\sigma}$ depending only on $\sigma$ such that
\begin{equation}\label{dist}
c_{\sigma} m\{\theta: N_{\sigma}u(\theta)>\lambda\}\leq P\{u^*>\lambda\}\leq C_{\sigma} m\{\theta: N_{\sigma}u(\theta)>\lambda\}
\end{equation}
for all $\lambda>0$.  Here, $m$ denotes the Lebesgue measure on the circle. 
\end{theorem}

From this distribution inequality and (\ref{BGS5}), it follows that for  $\Phi$ as in Theorem \ref{thm1}, 
\begin{equation}\label{BGS7}
\int_0^{2\pi} \Phi\left(N_{\sigma}u(\theta)\right) d\theta
\approx E\Phi(u^*)\approx E\Phi(v^*)\approx \int_0^{2\pi} \Phi\left(N_{\sigma}v(\theta)\right) d\theta
\end{equation} 
with constants depending only on $\sigma$ and $K$,
and  Theorem \ref{BGS2} follows. Here the middle $\approx$ is notation for (\ref{BGS5}) and the first and last $\approx$ follow from 
(\ref{BGS7}). 

These inequalities with $\Phi(x)=x^p$ for any $0<p<\infty$ give 
\begin{eqnarray*}
\sup_{0<r<1}\int_{0}^{2\pi} |F(re^{i\theta})|^p \,d\theta&\leq& 
2^p\big\{\int_0^{2\pi}|N_{\sigma}u(\theta)|^p \, d\theta +\int_0^{2\pi}|N_{\sigma}v(\theta)|^p \, d\theta\big\}\\
&\leq& C_{p, \sigma} \int_0^{2\pi}|N_{\sigma}u(\theta)|^p \, d\theta,
\end{eqnarray*}
and this proves the converse of (\ref{BGS1}). 

The paper \cite{Bur16} also proves versions of the above inequalities for the upper half-space $\bR^2_{+}$.  It is hard to overstate the influence of this paper, in both the probabilistic and analytic theory of Hardy spaces, in the 40 years since its publication. It would be virtually impossible to review here the literature that has its roots at least partially in \cite{Bur16}, and we do not make the attempt.  This paper was followed by the seminal paper of C. Fefferman and Stein \cite{FefSte},  and the theory of $H^p$ spaces (which had already been a hot area of research for many analysts) exploded from there.  For an account of $H^p$ theory and its connections to other subjects, see Stein \cite{Ste6}, Grafakos \cite{Gra1} and the literature referenced therein. 

\section{The good-$\lambda$ method}

We begin by showing how to prove Theorem \ref{thm1} using the good-$\lambda$ method. This involves proving distribution inequalities for the random variables $\tau^{1/2}$  and $B_{\tau}^*$ which, by a theorem that embodies the good-$\lambda$ method, immediately establish Theorem \ref{thm1}. We  first give the distribution inequalities  and then state the theorem which generates Theorem \ref{thm1} from them. The inequalities of the following theorem are from \cite{Bur20}. Although perhaps more ``elegant" and ``cleaner" than their original formulation in \cite{Bur14}, they may at first seem somewhat more mysterious.  

\begin{theorem}\label{thm2}
For all $0<\varepsilon<1$, $\delta>1$ and $\lambda>0$, 
\begin{equation}\label{good-one}
P\{B_{\tau}^*>\delta\lambda, \tau^{1/2}\leq\varepsilon\lambda\}\leq \frac{\varepsilon^2}{(\delta-1)^2}P\{B_{\tau}^*>\lambda\}
\end{equation}
and
\begin{equation}\label{good-two}
P\{\tau^{1/2}>\delta\lambda, B_{\tau}^*\leq\varepsilon\lambda\}\leq \frac{\varepsilon^2}{(\delta^2-1)}P\{\tau^{1/2}>\lambda\}.
\end{equation}
\end{theorem}

The following version of the good--$\lambda$ method from \cite{Bur20}, valid even for measure spaces of infinite measure, together with the inequalities just above, establishes Theorem \ref{thm1}.

\begin{theorem}\label{thm3}
Let $\Phi$ be as in the statement of Theorem \ref{thm1}. Suppose that
 $f$ and $g$ are nonnegative functions on a measure space  $(\Omega, \mathcal{A}, \mu)$, and $\delta>1$,
 $0<\varepsilon<1$, and $0<\gamma<1$  are real numbers such that
\begin{equation}\label{gengood}
\mu\{g>\delta\lambda, f\leq\varepsilon\lambda\} \leq\gamma \mu\{g>\lambda\},
\end{equation}
for every $\lambda>0$.
Let $\rho$ and $\nu$ be real numbers satisfying
$$
\Phi(\delta\lambda)\leq\rho\Phi(\lambda)\,\,\, \text{and}\,\, \ \Phi(\varepsilon^{-1} \lambda)\leq\nu\Phi(\lambda)$$
for every $\lambda>0$.
Finally, suppose that
\begin{equation}\label{assumption}
\rho\gamma<1\,\,\,  \text{and that}\,\,\, \int_{\Omega}\Phi(\min\{1, g\}) d\mu<\infty.
\end{equation}
Then
\begin{equation*}
\int_{\Omega}\Phi(g)d\mu \leq \frac{\rho\nu}{1-\rho\gamma} \int_{\Omega}\Phi(f)d\mu.
\end{equation*}
\end{theorem}

We note that in a space of finite measure, which is the case dealt with by Burkholder and Gundy in \cite{Bur14}, the assumption $\int_{\Omega}\Phi(\min\{1, g\}) d\mu<\infty$ is not needed, as it always holds.  For the general case it is known that this assumption is needed for the validity of the theorem;  see \cite{Jou} and \cite{MiyYab}.

Let us briefly outline the proof of this theorem following \cite{Bur20}. First note that 
\begin{eqnarray}
\mu\{g>\delta \lambda\} &\leq &\mu\{g>\delta\lambda, f\leq\varepsilon\lambda\} +\mu\{f>\epsilon\lambda\}\\
&\leq& \gamma \mu\{g>\lambda\}+\mu\{f>\epsilon\lambda\}\nonumber.
\end{eqnarray}
Multiplying by $\Phi'(\lambda)$ and integrating  in $\lambda$ gives
\begin{equation}
\int_{\Omega}\Phi\left(\frac{g}{\delta}\right)\,d\mu \leq \gamma\int_{\Omega}\Phi(g)\,d\mu+\int_{\Omega}\Phi\left(\frac{f}{\varepsilon}\right)\,d\mu.
\end{equation}
But
$$
\int_{\Omega}\Phi(g)\,d\mu =\int_{\Omega}\Phi\left(\delta\delta^{-1}g\right)\,
d\mu\leq \rho\int_{\Omega}\Phi\left(\frac{g}{\delta}\right)\,d\mu.
$$
This together with (\ref{assumption}) implies the theorem. (For full details, see 
\cite{Bur20}.)

The fact that both
$\frac{\varepsilon^2}{(\delta-1)^2}$ and $\frac{\varepsilon^2}{(\delta^2-1)}$
 in (\ref{good-one}) and (\ref{good-two}) go to zero as $\varepsilon$ goes to zero for any fixed $\delta$  is crucial in deriving Theorem \ref{thm1} from
 Theorem \ref{thm3}.  In the cases $\Phi(x)=x^p$, $0<p<\infty$, which  give the inequalities
 \begin{equation}
a_p \|\tau^{1/2}\|_p\leq \|B^{*}_{\tau}\|_p\leq A_p \|\tau^{1/2}\|_p,
 \end{equation}
the better the decay in  (\ref{good-one}) and (\ref{good-two}), as $\varepsilon$ goes to zero, the better the information obtained from Theorem \ref{thm2} on the constants $A_p$ and $a_p$ and the wider the applications to many other functionals involving $B^{*}_{\tau}$ and  $\tau$.    In \cite{Bur26},  Burkholder uses the good-$\lambda$ method to prove inequalities comparing $B^{*}_{\tau}$ and $\tau$ where $B$ is Brownian motion in $\bR^n$, $n\geq 1$; see also \cite{Dav3}. It follows from these argument that the optimal bounds on the right hand sides of  (\ref{good-one}) and (\ref{good-two}) are $C_1\exp\left(-\frac{(\delta-1)^2}{2\varepsilon^2}\right)$ and $C_2\exp\left(-C_3\frac{(\delta^2-1)}{\varepsilon^2}\right)$, respectively, where $C_1$, $C_2$ and $C_3$ are constants independent of $\delta$ and $\varepsilon$. An explicit value of $C_3$ can also be given.  (It should be noted that the proofs in \cite{Bur26} produce quantities which are denoted by  $R_n(\delta, \epsilon)$ and $L_n(\delta, \epsilon)$ in place of the exponentials.  An explicit computation for the case $n=1$ presented here provides the expressions given above; see \cite{BanMoo} for details.)  The proofs in \cite{Bur26} not only provide this ``gaussian" decay, which is best possible, but are extremely elegant and fairly simple, crystalizing, perhaps more than previous proofs, that what matters for the good-$\lambda$ method in the probabilistic setting is scaling and the strong Markov property.  The arguments in \cite{Bur26} can be extended to other functionals of Brownian motion
such as the maximal local time and good--$\lambda$ method proofs can be given of inequalities of
 Barlow and Yor \cite{BarYor1, BarYor2} which compare the norms of maximal local time to maximal functions and square functions, including various ratios of these quantities.  For more on this, see \cite{BanMoo}, \cite{Bass},  \cite{Dav2}, \cite{RevYor}, and references therein.

Let us also mention here that the paper \cite{Bur26} (see also \cite{Bur29}) contains several results which  have had a profound influence on probabilistic potential theory.  These include characterizations of Hardy spaces in domains of Euclidean space in terms of exit times of Brownian motion from the domains, and the computation of the  exact exponent $p$ for which the exit time $\tau$ for Brownian motion from an infinite cone in $\bR^n$ has finite $p$-moment. This exponent is given in terms of zeros of confluent hypergeometric functions.  This  led to many results on the exact growth (decay) of harmonic functions (and harmonic measure) on unbounded regions in $\bR^n$ and other applications to ``harmonic majorization" problems; see for example \cite{EssHalLewShe1}, \cite{EssHalLewShe2}.  There are also numerous other applications of these ideas to
conditioned Brownian motion and stable processes in more general cone-type regions.  The exact exponent of integrability of the exit time of these processes is related to the eigenvalues of the Laplace--Beltrami operator on the region of the sphere which generates the cone. This also leads to connections with heat kernels for ``singular" manifolds as in Cheeger \cite{Che1}. (See \cite{Deb1}, \cite{Deb2}, \cite{BanBog1}, \cite{BanSmi1},  and references therein.)

Motivated by their work on martingales in \cite{Bur14}, Burkholder and Gundy gave in \cite{Bur19} the first application of the  good-$\lambda$ method which did not involve probability in either the application of the method or the theorem proved. Denote the upper half-space by $\bR^{n+1}_+=\{(x, y): x\in \bR^n, y>0\}$ and the cone with vertex at $x$ and  aperture
$\alpha$ (for any $\alpha>0$) by $\Gamma_{\alpha}(x)=\{(\bar x, y)\in \bR^{n+1}_+ : |x-\bar x|<\alpha y\}$.   For any harmonic function $u$ defined on $\bR^{n+1}_+$, consider its nontangential maximal function (the upper half plane analogue of the function in (\ref{nontandisc}))
\begin{equation*}
N_{\alpha}u(x)=\sup_{( \bar{x}, y)\in \Gamma_{\alpha}(x)}\large|u(\bar{x}, y)\large|
\end{equation*}
and its Lusin area function
\begin{equation*}
A_{\alpha}u(x)=\left(\int_{\Gamma_{\alpha}(x)}y^{1-n}|\nabla u(\bar{x}, y)|^2 d\bar{x}\, dy\right)^{1/2}. 
\end{equation*}

These are harmonic function analogues of $Y^*$ and $\langle Y\rangle^{1/2}$ for martingales.  They have played a fundamental role in the development of harmonic analysis for the past seventy-five years. For an account of some of these applications;  see Stein \cite{Ste2, Ste3}.  Because of the  papers of Burkholder and Gundy \cite{Bur19} and Burkholder, Gundy and Silverstein \cite{Bur16} (already described above), the area and nontangental maximal functions and their use in  Littlewood--Paley theory have been inextricably connected to martingales and Brownian motion for the past 35-40 years.  In addition to the expository papers  of Burkholder in this volume on this subject, we refer the reader to Bass \cite{Bas}, Durrett \cite{Dur} and Gundy \cite{Gun1} and Varopoulos \cite{Var}.   The following version for the operators $A_{\alpha}$ and $N_{\alpha}$ of Theorem \ref{thm1} is proved in \cite{Bur19}. 

\begin{theorem} Let $\Phi$ be as in Theorem \ref{thm1}. There is a constant $C_1$ depending on $\alpha$, $n$ and $K$ such that
\begin{equation}
\int_{\bR^n} \Phi(A_{\alpha}u(x)) dx\leq C_1\int_{\bR^n} \Phi(N_{\alpha}(u)(x))dx.
\end{equation}
If the left hand side is finite, then $\lim_{y\to\infty}u(x, y)$ exists and does not depend on x. If $u$ is normalized so that  this limit is zero then there is a constant $C_2$ also depending on $\alpha$, $n$ and  $K$ such that
\begin{equation}
\int_{\bR^n} \Phi(N_{\alpha}u(x)) dx\leq C_2\int_{\bR^n} \Phi(A_{\alpha}(u)(x))dx.
\end{equation}
\end{theorem}

The good-$\lambda$ inequalities proved in \cite{Bur19} can be formulated as in Theorem \ref{thm2} as follows ($m$ denotes the Lebesgue measure in $\bR^n$).

\begin{theorem}  Let $0<\beta<\alpha$.  There is a constant $\gamma>1$ and  constants $C_1$ and $C_2$
depending only on $\alpha$, $\beta$ and $n$ such that for all  $0<\epsilon<1$ and $\lambda>0$,
\begin{equation}\label{goodharmonic-one}
m\{x\in \bR^n: N_{\beta}u(x)>\gamma\lambda, A_{\alpha}u(x)\leq\varepsilon\lambda\}\leq {C_1\,\varepsilon^2}\,m\{x\in \bR^n: N_{\beta}u(x)>\lambda\}
\end{equation}
and
\begin{equation}\label{goodharmonic-two}
m\{x\in \bR^n: A_{\beta}u(x)>\gamma\lambda, N_{\alpha}u(x)\leq\varepsilon\lambda\}\leq C_2\,\varepsilon^2
m\{x\in \bR^n: A_{\beta}u(x)>\lambda\}.
\end{equation}

\end{theorem}

Because of the ``local" nature of the proofs of these inequalities,  their proofs extends to measures other than Lebesgue measure which have densities that belong to the Muckenhoup $A_{\infty}$ class.  For this, the reader can see Gundy and Wheeden \cite{GunWhe} and Burkholder \cite{Bur32}.  

There are important applications of the good-$\lambda$ method to many other operators of fundamental importance in analysis and its applications such as the the Hardy--Littlewood maximal function and Calder\'on--Zygmund singular integrals (Coifman \cite{Coi} and Coifman-Fefferman \cite{CoiFef}),  to parabolic (heat equation) versions of $A_{\alpha}$ and $N_{\alpha}$ (Calder\'on-Torchisky \cite{CalTor}) as well as versions of these in the setting of Lipshitz domains and elliptic operators (Dahlberg \cite{Dal} and  Dahlberg-Jerison-Kenig \cite{DalJerKen}).  For more on  applications and further literature, see Garnett \cite{Gar}, Kenig \cite{Ken}, Stein \cite{Ste1} and Torchisky \cite{Tor}.

As in the case of martingales, the better the decay of the quantities on the right sides of (\ref{goodharmonic-one}) and (\ref{goodharmonic-two}) as $\varepsilon$ goes to zero, the wider the applications. The exponent $\varepsilon^2$ was improved to  $\varepsilon^k$ for any positive integer $k$ by R. Fefferman, Gundy, Silverstein and Stein \cite{FefGunSilSte} and by Murai and Uchiyma \cite{MurUch} to exponential $\exp{(-C/\varepsilon)}$ decay.  The full subgaussian  $\exp{(-C/\varepsilon^2)}$ decay as in the case of the martingale inequalities is proved in \cite{BanMoo}.

The original good-$\lambda$ inequalities in \cite{Bur19} imply the result of Privalov \cite{Pri}, Marcinkiewicz and Zygmund \cite{MarZyg}, Spencer \cite{Spe}, Calder\'on \cite{Cal1}, \cite{Cal2} and Stein \cite{Ste1} asserting that, except for sets of Lebesgue measure zero, the nontangential maximal functions and the Lusin area function are finite or infinite on the same sets and that on the sets where these are finite, the harmonic function has nontangential boundary limits; see \cite{Bur19}, \S3. In \cite{BanMoo}, the good-$\lambda$ gaussian decay inequalities are explored to prove a more quantitative version of these results which involve the { law of the iterated logarithm}.  The philosophy pioneered by Burkholder and Gundy, and Burkholder, Gundy and Silverstein, that these operators should be modeled after those for martingales with continuous paths and hence after stopped Brownian motion, drives the results presented in \cite{BanMoo}.

The paper \cite{Bur22} also studies boundary behavior of harmonic functions in the upper half space.  If $u$ is a harmonic function in the upper half-space  $\bR^{n+1}_+$ and $B$ is Brownian motion in $\bR^{n+1}_+$ and $\tau$ is the first time it leaves the upper half-space, then as before $u(B_{t\wedge\tau})$ is a martingale with Brownian maximal function  $u^*$ and Brownian square function $\langle u\rangle$.  In \cite{Bur22}, these quantities are used, in combination with techniques for (Doob's) conditioned Brownian motion, to prove that some (but not all) of the results of Privalov, Marcinkiewicz-Zygmund, Spencer, Calder\'on and Stein, can be obtained from the corresponding results for martingales. See \cite{Dur} for more on this.  The techniques in \cite{Bur22} have been used by other authors in different settings, see for example Brossard \cite{Bro1, Bro2} and Brossard and Chevelier \cite{BroCha1, BroCha2, BroCha3}.

\section{The Sharp Martingale Inequalities}

In the early eighties Burkholder turned his attention to sharp martingale inequalities and their extensions  to martingales taking values 
in Banach spaces.  The Banach space setting is discussed in Gilles Pisier's commentary in this volume.  In this section we discuss sharp inequalities and point out some of their implications in areas which are of current interest to many researchers.  In his seminal paper \cite{Bur39},  Burkholder proves the following theorem.

\begin{theorem}\label{sharpthm1} Let $f=\{f_n; n\geq 0\}$ be a martingale with difference sequence $d=\{d_n, n\geq 0\}$.  Suppose $1<p<\infty$ and let $p^*$  denote the maximum of $p$ and $q$ where $\frac{1}{p}+\frac{1}{q}=1$.  If $g$ is a  martingale transform of $f$ by a real predicable sequence $v$ uniformly bounded in absolute value by 1, then 
\begin{equation}\label{sharpeq1}
\|g\|_p\leq (p^*-1)\|f\|_p
\end{equation}
and the constant $(p^*-1)$ is best possible.   Furthermore, equality holds if and only if $p=2$ and $\sum_{k=0}^{\infty} v_k^2 d_k^2=\sum_{k=0}^{\infty} d_k^2$, almost surely.  
\end{theorem}

There are many other sharp martingale transform inequalities proved in \cite{Bur39}, including the following weak-type inequality. 
\begin{theorem}\label{sharpthm2} Let $1\leq p\leq 2$ and let $f$ and $g$ be as in Theorem \ref{sharpthm1}.  
Then 
\begin{equation}\label{sharpeq2}
\sup_{\lambda>0}\lambda^p\, P\{g^*>\lambda\}\leq \frac{2}{\Gamma(p+1)}\|f\|_p^p.
\end{equation}
The constant $ \frac{2}{\Gamma(p+1)}$ is best possible.  Furthermore, strict inequality holds if $0<\|f\|_p<\infty$ and $1<p<2$, but equality can hold if $p=1$ or 2. 
\end{theorem}

The case $p=1$ (with sharp constant 2) was proved in \cite{Bur31} and is not too difficult.  On the other hand, the proofs of Theorem \ref{sharpthm1} and  Theorem \ref{sharpthm2} (for the case when $1< p\leq 2$) are deep and difficult.  The proof of Theorem \ref{sharpthm1} (which after some preliminary work reduces to the case when the predicable sequence $\{v_k\}\in \{-1, 1\}$) rests on solving the nonlinear partial differential equation 
\begin{equation}
(p- 1)[yF_y - xF_x]F_{yy} - [(p - 1)F_{y} - xF_{xy}]^2 + x^2F_{xx}F_{yy} = 0
\end{equation}
for $F$ nonconstant and satisfying other conditions on a suitable domain of $\bR^2$.  Solving this equation leads to a system of five nonlinear differential inequalities with boundary conditions.  From this system, a function $u(x, y, t)$ is constructed in the domain 
$$
\Omega = \{(x, y, t)\in \bR^3: \big|\frac{x-y}{2}\big|^p <t\}
$$
with certain convexity properties from which, using the techniques of \cite{Bur35}, Burkholder proves that 
\begin{equation}
u(0,0,1)\|g_n\|_p^p\leq \|f_n\|_p^p
\end{equation}
for $1<p\leq 2$ and shows that $u(0,0,1)=(p-1)^p$.  This and duality give the bound $(p^*-1)$ in (\ref{sharpeq1}).  (The research announcement  in \cite{Bur37} contains a nice summary of the methods used in \cite{Bur39}.)

 Even today, the proofs in \cite{Bur39} seem extremely intricate. That Burkholder was able show that these PDE's have a solution with the important properties needed for the martingale inequalities is impressive.  A nice explanation of Burkholder's PDE and other ideas  in \cite{Bur39} in terms of the theory of Bellman functions was subsequently given by F. Nazarov, S. Treil and A. Volberg.  For this connection and some of the now very extensive literature on this subject, we refer the reader to  \cite{NazTre1, NazTreVol1, NazTreVol2, Vol1}.  Quoting from \cite{NazTreVol1}: ``It is really amazing that Burkholder was able to solve these PDEÕs: they are really complicated." 
  
 The Bellman function techniques have become a powerful tool to study sharp inequalities for many operators of great importance in harmonic analysis and its applications to PDE's and quasiconformal mappings.  In addition to the above already cited papers, see  \cite{Mel1, Mel2, NazVol1, SlaVas1, SlaVol1, SlaStoVas1, PetVol1, VasVol1}.

In a series of papers following \cite{Bur39}, which included many applications to various other sharp inequalities for discreet  martingales and stochastic integrals, Burkholder simplified the proofs in \cite{Bur39} considerably by giving  explicit expressions for his ``magical" functions.  In the Bellman function language of Nazarov and Volberg, Burkholder gives an explicit expression for the ``true" Bellman function of the above PDE.  Quoting from \cite{NazTre1}, ``the most amazing thing is that the true Bellman function is known! This fantastic achievement belongs to Burkholder."  Explicit solutions to Bellman problems that arise in many of the applications to harmonic analysis are often extremely challenging problems.  For more on this, see \cite{Mel1}, \cite{Mel2}, \cite{VasVol1}, and especially the recent paper \cite{VasVol2} which contains a nice treatment, based on Monge-Amp\`ere equation, on how to solve many Bellman equations, including Burkholder's. 

We now recall the following generalization of Theorem \ref{sharpthm1}, which is proved by Burkholder in \cite{Bur45} using the explicit form of his  function $U$. 

 \begin{theorem}\label{sharpthm3} Let $\bH$ be a (real or complex) Hilbert space. For $x\in \bH$, let $|x|$ denote its norm.  Let $f=\{f_n\}_{n=0}^{\infty}$ and $g=\{g_n\}_{n=0}^{\infty}$ be two $\bH$-valued martingales on the same filtration with martingale difference sequences $d$ and $e$, respectively, and satisfying  $|e_k|\leq |d_k|$ pointwise for all $k\geq 0$.  Then, with $p$ and $p^*$ as in Theorem \ref{sharpthm1}, 
$\|g\|_p\leq (p^*-1)\|f\|_p$. The constant $(p^*-1)$ is best possible and equality holds (in the case $0<\|f\|_p<\infty$) if and only if $p=2$ and $|e_k|=|d_k|$ almost surely, for all $k\geq 0$. 
\end{theorem}

To prove this inequality, Burkholder considers the function $V:\bH\times\bH\to \bR$\,    defined by 
\begin{equation}\label{v}
V(z, w)=|w|^p-(p^*-1)^p|z|^p.
\end{equation}
The goal is then to show that $EV(f_n, g_n)\leq 0$.  Burkholder then introduces the  function 
\begin{equation}\label{u}
U(z, w)=\alpha_p\left(|w|-(p^*-1)|z|\right)\left(|z|+|w|\right)^{p-1}
\end{equation}
where 
$$
\alpha_p=p\left(1-\frac{1}{p^*}\right)^{p-1}
$$
and proves that this function satisfies the following properties: 
\begin{eqnarray}\label{sub1}
V(z, w)&\leq &U(z, w)\,\,\,  \text{for all}\,\,\, w, z\in \bH,\label{sub2}\\
EU(f_n, g_n) &\leq& EU(f_{n-1}, g_{n-1}), \,\,\, n\geq 1,\label{sub3}\\
EU(f_0, g_0)&\leq& 0. 
\end{eqnarray}

There have been many applications of these ideas to martingales, harmonic functions (including differential subordination of harmonic functions) and singular integrals. Many of these results are due to Burkholder himself; see for example his work in \cite{Bur40, Bur42, Bur43, Bur44, Bur46, Bur47, Bur48, Bur50, Bur51, Bur52, Bur54}.  Some other applications (including many recent ones) are contained in \cite{Cho3, Cho4, Cho2, Cho1, Ham2, Ham1, Jan1, Osc1, Osc2, Osc3, Osc4, Osc5, Osc6, Osc7, Suh1, Wan2, Wan3, Wan4, Wan1}.   

\subsection{The $(p^*-1)$ constant in analysis}

No sooner had Burkholder's paper \cite{Bur39} appeared identifying the  $L^p$ norm of martingale transforms as $(p^*-1)$, than the connection (at least superficially at first) to a conjecture of T. Iwaniec \cite{Iwa1} concerning the $L^p$ norm of the Beurling--Ahlfors operator was noticed by several researchers.   The Beurling--Ahlfors operator is a singular integral operator (Fourier multiplier) 
on the complex plane $\bC$ (or $\bR^2$) defined on $L^p(\bC)
\cap L^2(\bC)$, $1<p<\infty$, by 
\begin{equation}\label{BAdefinition}
Bf(z) = -\frac{1}{\pi}\textrm{p.v.}\int_{\bC}\frac{f(w)}{(z-w)^2}dm(w),  \hskip.3cm \widehat{Bf}(\xi)=\frac{\overline{\xi}^2}{|\xi|^2}\hat f(\xi), \hskip.3cm \xi\not=0\in \bC.
\end{equation}

This operator (which incidentally can be written in terms of second order  Riesz transforms in the plane as $B=R_2^2-R_1^2+2iR_1R_2$) is of fundamental importance in several areas of mathematics including PDE and 
 the geometry of quasiconformal mappings \cite{ Ast1, AstIwaMar1, DonSul1, Iwa1, Iwa2, Iwa3, IwaMar1, PetVol1}.  As a Calder\'on--Zygmund singular integral, it is bounded on $L^p(\bC)$, for $1<p<\infty$. The computation of its operator norm, $\|B\|_p$, has been an open problem for almost thirty years. 
In  \cite{Leh1}, Lehto showed that $\|B\|_p \geq  p^*-1$. Inspired in part by the celebrated Gehring--Reich
 conjecture \cite{GerRei} on the area distortion of qasiconformal mappings in the plane 
 (proved by K. Astala  \cite{Ast1}), T. Iwaniec conjectured in \cite{Iwa1} that $\|B\|_p = p^*-1$.  For some of the connections to quasiconformal mappings, see K. Astala, T. Iwanienc and G. Martin  \cite{AstIwaMar1}.
  
 Although the $L^p$--boundedness of very general singular integrals and Fourier multiplier operators can be proved using martingale transforms (see for example Burkholder \cite{Bur38} and McConnell \cite{McC1}), obtaining precise information on their norms requires more exact representations of the operators in terms of martingales.  In the groundbreaking paper \cite{GunVor1}, Gundy and Varopoulos gave a representation for Riesz transforms in terms of stochastic integrals arising from composing the harmonic extension of the function $f$ with Brownian motion in the upper half-space (as in \S3 above).  The Burkholder results in \cite{Bur39} and the Gundy--Varopoulos representation  provided hope that the Iwaniec conjecture could be derived from Burkholder's theorem \cite{Bur39}. Unfortunately, to implement this approach certain stochastic integral versions of Burkholder's inequalities were needed and these did not follow, in any direct way, from the results in \cite{Bur39}. Hence, those interested in this approach had to wait. The wait was over when Burkholder gave the explicit expression for his function $U$, which together with the It\^o formula leads to the desired sharp stochastic integral inequalities that arise when applying the Gundy--Varopoulos formula to the Beurling-Ahlfors operator; see \cite{BanWan3}, \cite{BanWan2}, \cite{BanWan1}, \cite{Wan1}.  This approach was employed in \cite{BanWan3} to show that $\|Bf\|_p\leq 4(p^*-1)\|f\|_p$.  While not the desired bound, this paper provided the first explicit bound for $\|B\|_p$.  In \cite{NazVol2}, Nazarov and Volberg improved this bound to $2(p^*-1)$ by proving a Littlewood-Paley inequality using Bellman functions.  However, the construction of their Bellman function (for which no explicit expression is known as of now) depends on Burkholder's sharp inequality for Haar martingales. In \cite{BanMen1}, the proof in \cite{BanWan3} is redone using the heat kernel in place of the Poisson kernel to obtain the $2(p^*-1)$ bound.   Both the stochastic integral techniques in \cite{BanMen1} and   \cite{BanWan3}, and the Nazarov-Volberg Bellman approach \cite{NazVol2}, have had many other applications, including applications to bounds for the Beurling-Ahlfors operator in several dimensions (first studied in  \cite{DonSul1} and  \cite{IwaMar0}) and Riesz transforms for the Ornstein-Uhlenbeck process.  Some of these applications are contained in  \cite{Arc1, ArcLi1, BanBog2,  BanJan1, BanLin1,  BanMen1,   DraPetVol1, DraVol1, DraVol2, DraVol3, DraVol4, GeiSmiSak1, Hyt}. 
 The best bound for $\|B\|_p$ as of now, proved in \cite{BanJan1}, is $1.575(p^*-1)$. This bound is obtained by an improvement of Burkholder's $(p^*-1)$ inequality for complex valued martingales which have some additional orthogonality structure between its real and imaginary parts. The explicit expression of the  Burkholder function $U$ is crucial for this proof.

 In \cite{BanWan3} and \cite{BanWan2} the techniques of Burkholder are adapted to show that for martingales under differential subordination and orthogonality, the Burkholder constant $(p^*-1)$ can be replaced by the Pichorides \cite{Pic1} Hilbert transform constant $\cot({\frac{\pi}{2p^*}})$.  From this and the Gundy--Varopoulos stochastic representation, it follows that for the Riesz transforms $\|R_j\|_p\leq \cot({\frac{\pi}{2p^*}})$, $j=1,2,\dots n$.  This result, and other good estimates of the $L^p$-norms of singular integral operators with odd kernels, can also be obtained by the classical method of rotations; see \cite{DuoRub1}, \cite{IwaMar1},  \cite{Pis1}, and especially \cite{IwaMar1} where it is shown that in fact $\|R_j\|_p=\cot({\frac{\pi}{2p^*}})$. It is also interesting to note here that the Beurling-Ahlfors operator plays the role of the Hilbert transform in the so-called {complex method of rotation} and this can also be use to obtain estimates of the $L^p$-norms of certain singular integral operators with even kernels in terms of $\|B\|_p$. 
We refer the reader to \ \cite{IwaMar1} where this technique was introduced;  see also, \cite{IwaMar2}.

 In the recent paper \cite{GeiSmiSak1}, Geiss, Mongomery-Smith and  Saksman combined the estimates in \cite{BanMen1, NazVol2} and the arguments used by Bourgain in \cite{Bou1} which showed that Burkholder's UMD property is  equivalent to the boundedness of the Hilbert transform, to prove that the $L^p(\bR^n)$ operator norm of  $2R_jR_k$, $j\not= k$, is $(p^*-1)$.  That is, $\|2R_jR_k\|_p=(p^*-1)$, $j\not= k$.  This beautiful and surprising result gives the first example of a singular integral whose $L^p$ norm is exactly that of martingale transforms.  
In particular, in the plane the two ``components" of the Beurling--Ahlfors operator $B$ satisfy $\|R_2^2-R_1^2\|_p=(p^*-1)$ and $\|2R_1R_2\|_p=(p^*-1)$.   The proof in \cite{GeiSmiSak1} adapts to other combinations of $R_j$ and $R_k$; see for example \cite{VasVol2}. 

 \subsection{Rank-one convexity and quasiconvexity}

 It is well known (\cite{AstIwaMar1}) that proving $\|Bf\|_p\leq (p^*-1)\|f\|_p$ for all $f\in L^p(\bC)$, $1<p<\infty$, is equivalent to proving 
 \begin{equation}\label{BA-2}
 \|\partial f\|_p\leq (p^*-1)\|\overline{\partial}f\|_p, \,\,\,\, 1<p<\infty, 
 \end{equation}
  for all smooth functions $f$ of compact support ($f\in C_0^{\infty}(\bC)$) where (with $z=x+iy$)
\begin{equation}\label{cauchy}
\partial=\left(\frac{\partial f}{\partial x}-i\frac{\partial f}{\partial y}\right) \,\,\, \text{and}\,\,\, \overline\partial=\left(\frac{\partial f}{\partial x}+i\frac{\partial f}{\partial y}\right)
\end{equation}
 are the Cauchy-Riemann operators in the complex plane.  Viewed in terms of the function $V$ in (\ref{v}), (\ref{BA-2}) is the same as proving that
 \begin{equation}\label{BA-V} 
 \int_{\bC}V(\overline\partial f, \partial f)dm(z)\leq 0,\,\,\,\,  f\in C_0^{\infty}(\bC). 
 \end{equation}
 Since by Burkholder (\ref{sub1}), $V(z, w)\leq U(z, w)$ for all $w, z\in \bC$, it is natural to conjecture that 
 \begin{equation}\label{BA-U} 
 \int_{\bC}U(\overline\partial f, \partial f)dm(z)\leq 0,\,\,\,\,  f\in C_0^{\infty}(\bC). 
 \end{equation}
 This conjecture, which arose from the work in \cite{BanWan3}, is written down as a question ({\it Question 1}) in \cite{BanLin1}.  
 The conjectured inequality (\ref{BA-U}) and the convexity properties (listed below) satisfied by the function $U$ lead to another unexpected connection and application of Burkholder's powerful ideas.   
 
 Denote  by $\M^{n\times m}$ the set of all $n\times m$ matrices with real entries.  
 The function $\Psi:\M^{n\times m}\to \bR$  is said to be  {rank-one convex} if for each
$A, B \in \M^{n\times m}$ with rank $B = 1$, the function
\begin{equation}\label{rank-one}
h(t)=\Psi(A+tB), \, \, \,\, t\in \bR
\end{equation}
is convex.  The function is said to be {quasiconvex} if  it is locally integrable and for each $A\in \M^{n\times m}$, 
bounded domain $\Omega\subset \bR^n$ and each compactly supported Lipschitz function $f : \Omega\to \bR^m$ we have 
  \begin{equation}\label{quasiconvex}
  \Psi(A)\leq \frac{1}{|\Omega|}\int_{\Omega} \Psi\left(A+Df(x)\right)\,dx. 
  \end{equation}
 where $Df$ is the Jacobian matrix of $f$.  
 
These properties arise in many problems in the calculus of variations, 
especially in efforts to extend the so called ``direct method" techniques from convex energy functionals to nonconvex.  They were introduced by C.B. Morrey (see \cite{Mor2}) and further developed by J. Ball \cite{Bal1}.  For more (much more) on the relationship
 between these properties and their consequences in the direct method of the calculus of variations, we refer to \cite{Dac1}.  If $n = 1$ or $m = 1$, then $\Psi$ is quasiconvex or rank one convex if and only if it is convex. If $m\geq 2$ and 
and $n\geq 2$, then convexity $\Rightarrow$ quasiconvexity $\Rightarrow$ rank-one convexity. (See \cite{Dac1} where the notion of polyconvexity  which lies ``between" convexity and quasiconvexity is also discussed.)  
 In 1952, Morrey \cite{Mor1} conjectured that rank-one convexity does not imply quasiconvexity when
both $m$ and $n$ are at least 2. In 1992,  \v{S}ver\'ak \cite{Sve2} proved that Morrey's conjecture is correct if $m\geq 3$  and $n\geq 2$. The cases $m=2$ and $n\geq 2$ remains open.  One of the difficulties with these notions of convexity 
 is that it is in general very difficult to construct nontrivial, interesting, examples of such functions.

 Enter Burkholder's function $U$. It is proved in \cite{Bur45} that for all $z,\ w,\ h,\ k\in \bC$ with  $|k|\leq |h|$, the function  
$t\to U (z+th,\ w+tk)$  is concave in $\bR$ , or equivalently that $t\to -U (z+th,\ w+tk)$ is convex in $\bR$ . The concavity property of $t\to U (z+th,\ w+tk)$ is crucial in the proof of the properties in (\ref{sub1})--(\ref{sub3}). Properly interpreted, this convexity property of $U$ is equivalent to rank-one convexity. Define the function 
$
\Gamma\colon \M^{2\times 2}\to \bC\times \bC$
by 
$$
\Gamma\left(\begin{array}{cc}a & b \\c & d\end{array}\right)=(z, w),
$$ where 
$z=(a+d)+i(c-b)$ and  $w= (a-d)+i(c+b)$ and set 
$\Psi_{U}= -U\circ \Gamma$.  It follows easily from the convexity property of $t\to -U (z+th,\ w+tk)$, for $z,\ w,\ h,\ k\in \bC$ with  $|k|\leq |h|$, that the function $\Psi_U$ is rank-one convex. (See \cite{BanLin1} for full details.)
Now, if $f=u+iv \in C_{0}^{\infty}(\bC)$, then 
\begin{equation}
Df=\begin{pmatrix}
     u_x &  u_y  \\
     v_x &  v_y
\end{pmatrix}
\end{equation}
and 
\begin{equation}
\Psi_{U} \left(Df\right)
=-U\left({\overline \partial f},\ {\partial f}  \right).
\end{equation}
Thus quasiconvexity of $\Psi_{U}$ (at $0\in \bR^{2\times 2}$) is equivalent to 
\begin{equation}\label{quasi} 
 0\leq -\int_{\text{supp f}}U(\overline\partial f, \partial f)dm(z),
 \end{equation}
which is equivalent to (\ref{BA-U}).  Thus the question:  Is Burkholder's  function $U$ also quasiconvex in the sense that $\Psi_U$ is quasiconvex?  If the answer is ``yes", then the Iwaniec 1982 conjecture follows.  If the answer is ``no" then the Morrey 1952 conjecture follows for the important case $n=m=2$.   

The article by A. Baernstein and Montgomery-Smith \cite{BaeSmt1} presents various connections between the function $U$ and another function $L$ used by Burkholder to prove sharp weak--type inequalities for martingales and harmonic functions under the assumption of differential subordination, \cite[p. 20]{Bur46}.  This function $L$ was subsequently independently rediscovered by  \v{S}ver\'ak in \cite{Sve3}.  For more on these connections, we refer the reader to  \cite[pp. 518-523]{AstIwaMar1}, \cite{Iwa3}, \cite{Sve1, Sve2, Sve3} and \cite{VasVol2}. 

Finally, we give a brief account of recent developments in which quasiconformal mappings (also nonlinear hyperelasticy) and  Burkholder's theory on sharp martingale inequalities share common problems of compelling interest. (We refer the reader to \cite{Iwa3} and \cite{AstIwaMar1} for details.) By definition, a weakly differentiable mapping $f :\Omega \rightarrow \bR^n$ in a domain  $\Omega\subset \bR^n$ (also referred to as hyperelastic deformation) is said to be $K$-quasiregular, $1\leq K < \infty ,$ if its Jacobian matrix  $D\!f(x) \in \M^{n\times n}$ (deformation gradient) satisfies the distortion inequality
\begin{equation}\label{distor}
 |D\!f(x)|^n  \leqslant K\, \det D\!f(x),\;\quad \textrm{where}\;\;\;|D\!f(x)|  = \max_{|v| =1} \; |D\!f(x) v|. 
\end{equation}
 
 The $L^p$-integrability of the derivatives of $K$-quasiregular mappings relies on a general inequality which is opposite to the distortion inequality in an average sense. More precisely, 
\begin{equation} \label{1}
 \int_{\mathbb R^n} \left \{\,|D\!F(x)|^n \;-\; K\, \det D\!F(x) \,\right \}\cdot |DF(x)|^{p-n}\; \textrm d x \;\geqslant  \; 0, 
\end{equation}
for all mappings $F \in \W^{1,p}(\bR^n, \bR^n)$ with the Sobolev exponents $p$ in a certain interval $ \alpha(n,K)  < p < \beta(n,K)$, where $ \alpha(n,K)  < n <\beta(n,K)$. 
Iwaniec (\cite{Iwa3}) conjectured that the largest such interval is: 
\begin{equation}\label{2}
\alpha (n,K)  = \frac{nK}{K+1} < p < \frac{nK}{K-1} = \beta(n,K).
\end{equation}
Iwaniec (see again \cite[pp. 518-523]{AstIwaMar1} and \cite{Iwa3}) then observed that in dimension $n=2\,$ the integrand in  (\ref{1}) is none other than the Burkholder's function $U$ (modulo constant factor), thus rank-one convex for all exponents $p$ in (\ref{2}). Inspired by Burkholder's results he proved, in every dimension $n\geqslant 2$, that (\ref{2}) defines precisely the range of the exponents $p$  for which the integrand in (\ref{1}) is rank-one convex; see \cite{Iwa3}. Now, it may very well be that Iwaniec's $n$-dimensional analogue of Burkholder's integral is also quasiconvex and, conjecturally, that (\ref{1}) holds for all $p$ in the range (\ref{2}). This would give a completion of the  $L^p$-theory of quasiregular mappings in space.

While it is not clear at this point that martingale techniques will
produce the sharp bound $(p^*-1)$ for $\|B\|_p$, it seems likely that
the solution to the Iwaniec conjecture will somehow involve the
Burkholder function $U$. Also, in higher dimensions it is plausible
that Burkholder's vision and his sharp martingale inequalities will
contribute the  creation of a viable $L^p$-theory of quasiregular
mappings with far reaching applications to geometric function theory in
$\bR^n$ and, in particular, mathematical models of nonlinear
hyperelasticity.  What is certainly clear is that as of now all approaches (stochastic integrals and
Bellman functions) which have produced concrete bounds close to the
conjectured bound for $\|B\|_p$ rest on the fundamental ideas of
Burkholder originally conceived to prove sharp martingale inequalities. These  ideas have led to deep and
surprising connections in areas of analysis and PDE's where this and
other singular integrals operators and maximal functions play an
important role and which on the surface seem far removed from martingales. 
More than twenty five years after the publication of
\cite{Bur39}, the techniques and ideas in this paper are still being
explored by many mathematicians in different fields. This is indeed a
landmark paper, one of many in Burkholder's list of publications.

\end{document}